\documentclass[11pt,leqno]{amsart}
\usepackage{amsmath,amsfonts,amssymb,amsthm}

\theoremstyle{plain}
\newtheorem{theorem}{Theorem}

\newtheorem{lemma}[theorem]{Lemma}

\theoremstyle{definition}

\newcommand{\barint}{
\rule[.036in]{.12in}{.009in}\kern-.16in \displaystyle\int }

\newcommand{\barcal}{\mbox{$ \rule[.036in]{.11in}{.007in}\kern-.128in\int $}}

\newcommand{\B}{\mathbb B}
\newcommand{\R}{\mathbb R}

\newcommand{\Z}{\mathbb Z}

\newcommand{\Heis}{\mathbb H}
\newcommand{\Sph}{\mathbb S}

\newcommand{\lip}{{\rm Lip}\, }

\newcommand{\cH}{\mathcal H}

\makeatletter
\def\Ddots{\mathinner{\mkern1mu\raise\p@
\vbox{\kern7\p@\hbox{.}}\mkern2mu
\raise4\p@\hbox{.}\mkern2mu\raise7\p@\hbox{.}\mkern1mu}}
\makeatother

\numberwithin{theorem}{section} \numberwithin{equation}{section}

\title[Lipschitz homotopy groups]{The ($n+1$)-Lipschitz homotopy group of the Heisenberg group $\Heis^n$}

\author{Piotr Haj{\l}asz}

\address{P.\ Haj{\l}asz: Department of Mathematics, University of Pittsburgh, 301
  Thackeray Hall, Pittsburgh, PA 15260, USA, {\tt hajlasz@pitt.edu}}

\thanks{P.H.\ was supported by NSF grant DMS-1500647.}

\begin{document}

\sloppy


\begin{abstract}
We prove that for $n\geq 2$, the Lipschitz homotopy group $\pi_{n+1}^{\lip}(\Heis^n)\neq 0$ of the Heisenberg group $\Heis^n$
is nontrivial.
\end{abstract}

\maketitle

\section{Introduction}

The Lipschitz homotopy groups of a metric space $\pi_n^{\lip}(X)$ are defined in the same way as the 
classical homotopy groups, with the exception that both
the maps and homotopies are required to be Lipschitz. We emphasize that we make no restriction on the Lipschitz constants. In particular,
we do not require the Lipschitz constant of a homotopy to be comparable to the Lipschitz constant of maps that are Lipschitz homotopic.
The notion of Lipschitz homotopy groups was introduced in \cite{DHLT} with the main purpose of studying the Lipschitz homotopy groups of the
Heisenberg group (equipped with the Carnot-Carath\'eodory metric). 
Note that the classical homotopy groups of the Heisenberg group are zero since the space is homeomorphic to a Euclidean space. The following results are known
$$
\pi_k^{\lip}(\Heis^n)
\left\{
\begin{array}{cccc}
=0      & \mbox{if} & 1\leq k<n          & \text{\cite{WY1}},\\
\neq 0  & \mbox{if} & k=n                & \text{\cite{BF}},\\
=0      & \mbox{if} & n=1, k\geq 2       & \text{\cite{WY2}},\\
\neq 0  & \mbox{if} & n=2\ell, k=4\ell-1 & \text{\cite{HST}}.\\
\end{array}
\right.
$$
The main result of this paper adds one more case to this list.
\begin{theorem}
\label{main}
For $n\geq 2$,
$\pi_{n+1}^{\lip}(\Heis^n)\neq 0$.
\end{theorem}
The case $n=2$ is already contained in \cite{HST}, but with a very different proof.
As we will see, the result has a very short proof and it is somewhat surprising that it has not been noticed before.

The proofs that $\pi_n^{\lip}(\Heis^n)\neq 0$ and $\pi_{4\ell-1}^{\lip}(\Heis^{2\ell})\neq 0$ employ the fact that the groups 
$\pi_n(\Sph^n)$ and $\pi_{4\ell-1}(\Sph^{2\ell})$ are infinite so one can use differential forms to detect non-trivial elements
in these groups (using degree and the Hopf invariant respectively). To the best of my knowledge, Theorem~\ref{main} for $n\geq 3$,
is a first result showing non-trivial Lipschitz homotopy groups of the Heisenberg groups when the corresponding 
homotopy groups of the spheres
are finite, $\pi_{n+1}(\Sph^n)=\Z_2$ for $n\geq 3$.

When $k>1$, Lipschitz homotopy groups are abelian, but they are usually uncountably generated; for a detailed proof that
$\pi_n^{\lip}(\Heis^n)$ is uncountably generated, see \cite{DHLT}. However, very little is known about the structure of these groups.
For example it is not known if $\pi_n^{\lip}(\Heis^n)$ has elements of finite order. Nevertheless, the fact that a Lipschitz 
homotopy group is nontrivial provides a lot of information about the structure of Lipschitz maps into the Heisenberg group.

In order to show that $\pi_k^\lip(\Heis^n)\neq 0$, we need to construct a Lipschitz mapping $f:\Sph^k\to\Heis^n$ that has no Lipschitz extension
$F:\B^{k+1}\to\Heis^n$, see comments following Definition 4.1 in \cite{DHLT}. So far, the only known method of proving that $\pi_k^{\lip}(\Heis^n)\neq 0$ is based on the idea
described below.

The following fact is well known, see \cite{BF,DHLT,EES}.
\begin{lemma}
\label{L1}
For $n\geq 1$, there is a smooth and horizontal embedding of the sphere $\Phi:\Sph^n\to\Heis^n$,
that is $\Phi:\Sph^n\to\R^{2n+1}$ is a smooth embedding and the tangent space to the embedded
sphere $\Phi(\Sph^n)$ is contained in the horizontal distribution of $\Heis^n$.
It follows that the embedding $\Phi:\Sph^n\to\Heis^n$ is bi-Lipschitz.
\end{lemma}

Now if $\pi_k(\Sph^n)\neq 0$ and $f\in C^\infty(\Sph^k,\Sph^n)$, $0\neq [f]\in \pi_k(\Sph^n)$, then, in some cases, one can prove that
if $\Phi$ is any map that satisfies the conclusion of Lemma~\ref{L1}, then
$\Phi\circ f:\Sph^k\to\Heis^n$ does not admit a Lipschitz extension $F:\B^{k+1}\to\Heis^n$. This was the method used in all of the
known cases of $\pi_k^{\lip}(\Heis^n)\neq 0$ and it is also a method used in this paper. For this reason I conjectured that if
$\pi_k(\Sph^n)\neq 0$, then $\pi_k^\lip(\Heis^n)\neq 0$. While the conjecture is still open, in many cases the approach 
to a proof would have to be very different from the one described here, because of the following result of Wenger and Young \cite[Theorem~1]{WY2}.
\begin{theorem}
\label{WY}
If $f:\Sph^k\to\Sph^n$ is Lipschitz and $n+2\leq k<2n-1$, then $\Phi\circ f:\Sph^k\to\Heis^n$ can be extended to a Lipschitz map
$F:\B^{k+1}\to\Heis^n$.
\end{theorem}
The result is sharp: it does not extend to the case of $k=2n-1$ when $n$ is even. Indeed,
if $n$ is even and $f:\Sph^{2n-1}\to\Sph^n$ has a non-zero Hopf invariant, then as was shown in \cite{HST}, the map
$\Phi\circ f:\Sph^{2n-1}\to\Heis^n$ does not admit a Lipschitz extension $F:\B^{2n}\to\Heis^n$ and hence $\pi_{2n-1}^{\lip}(\Heis^n)\neq 0$. 
Also the proof of Theorem~\ref{main} shows that Theorem~\ref{WY} cannot be extended to the case $k=n+1$ which is somewhat surprising in view of \cite[Theorem~2]{WY2}.

Note that Theorem~\ref{WY} shows that the method described above of proving that $\pi_k^\lip(\Heis^n)\neq 0$ cannot be used in the
range $n+2\leq k<2n-1$.

\section{Proof of Theorem~\ref{main}}

Let $\Phi:\Sph^n\to\Heis^n$ be a smooth and horizontal embedding of the sphere into the Heisenberg group as in Lemma~\ref{L1}.
Recall that
$$
\pi_{n+1}(\Sph^n)=
\left\{
\begin{array}{ccc}
\Z    & \mbox{if $n=2$,}\\
\Z_2
&   \mbox{if $n\geq 3$.}
\end{array}
\right.
$$
Let $f\in C^\infty(\Sph^{n+1},\Sph^n)$ be a smooth map such that $0\neq [f]\in\pi_{n+1}(\Sph^n)$.
We will show that the map
$$
g=\Phi\circ f:\Sph^{n+1}\to\Heis^n
$$
does not admit a Lipschitz extension 
$$
G:\B^{n+2}\to\Heis^n
$$
which will prove that $0\neq [g]\in \pi_{n+1}^\lip(\Heis^n)$.
Suppose to the contrary that $G:\B^{n+2}\to\Heis^n$ is a Lipschitz extension of $g$. Since the Heisenberg group
is purely unrectifiable in dimensions $k\geq n+1$, \cite{ambrosiok,magnani} (see \cite{BHW} for a different proof),
we conclude that
$\cH^{n+2}(G(\B^{n+2}))=0$, where $\cH^{n+2}$ stands for the Hausdorff dimension with respect to the Carnot-Carath\'eodory
metric in $\Heis^n$. To arrive to a contradiction it suffices to show that $\cH^{n+2}(G(\B^{n+2}))>0$.

We will need the following result of Gromov \cite[Theorem~ 3.1.A]{gromov}, see also \cite[Lemma~22 and Theorem~6]{pansu}.
\begin{lemma}
If a set $E\subset\Heis^n$ has topological dimension at least $n+1$, then $\cH^{n+2}(E)>0$.
\end{lemma}
In fact the statement of Gromov's result only says that the Hausdorff dimension of $E$ is greater than or equal $n+2$
so one cannot conclude from the statement of Gromov's result that $\cH^{n+2}(E)>0$.
However, the proof of Gromov's result clearly shows that
$\cH^{n+2}(E)>0$, see the proof of Lemma~22 in \cite{pansu}.

Thus it suffices to show that the topological dimension of $G(\B^{n+2})$ is greater than or equal to $n+1$.
Suppose to the contrary that $\dim_{\rm top} G(\B^{n+2})\leq n$. We will need the following classical result
\cite[Theorem~VI.4]{HW}.
\begin{lemma}
A separable metric space $X$ has topological dimension less than or equal $n$ if and only if 
for each closed set $C\subset X$ and a continuous map $h:C\to \Sph^n$, there is a continuous extension
$H:X\to\Sph^n$ of $h$.
\end{lemma}
Let $h=\Phi^{-1}:\Phi(\Sph^{n})\to\Sph^n$. Since $\Phi(\Sph^n)$ is a compact subset of $G(\B^{n+2})$ and 
$\dim_{\rm top} G(\B^{n+2})\leq n$, there is a continuous extension 
$H:G(\B^{n+2})\to\Sph^n$ of $h$. Note that
$H\circ G:\B^{n+2}\to\Sph^n$ is a continuous extension of
$$
H\circ G\Big|_{\Sph^{n+1}}=f:\Sph^{n+1}\to\Sph^n
$$
which contradicts the assumption that $0\neq [f]\in\pi_{n+1}(\Sph^n)$.
The proof is complete.
\hfill $\Box$

\end{document}